\input amstex
 \input amsppt.sty
 \documentstyle{amsppt}
 \NoRunningHeads
 \NoBlackBoxes
 \topmatter
 \title Irreducible Components in an Algebraic Variety of Representations
 of a family of  One-relator Groups
 \endtitle
 \author Sal Liriano
 \endauthor
 \abstract
 Given a finitely generated  group  $G$, the set  $Hom(G, SL_2 \bold C)$
 inherits  the
 structure  of  an  algebraic  variety   $R(G)$  called
 the {\it  representation variety \/}
 of $G$. This algebraic variety is an invariant of $G$.
 Let  $$G_{pt}=\langle a, b; a^p= b^t \rangle, $$
  where $p$, $t$ are integers greater than one.
 In this paper a formula is produced yielding the number of four dimensional
 irreducible components of the  affine  algebraic  variety
 $R(G_{pt})$. A direct consequence of the main theorem of this paper is that if $K$ is a  torus
 knot,
 then  its genus   equals   the   number   of   four   dimensional
 components  of the representation variety, $R(G_{pt})$, of its corresponding knot group $G_{pt}$.
 \endabstract
 \endtopmatter
 \document
 \head Introduction
 \endhead

  Let $G$ be a  finitely generated  group generated by
 $X = \{ x_1,x_2,\dots ,x_n \}$;
 then the set $Hom(G, SL_2 \bold C )$ can  be
 endowed  with  the  structure
 of
 an  affine algebraic variety
 (see \cite {BG}, \cite {LM})
 here denoted by $R(G)$. The algebraic variety $R(G)$  is
 independent of the
 choice of finite set $X$ of generators for $G$. In other  words,
 if a second finite  set  of  generators  $X' $  for  $G$
 is chosen, the
 resulting algebraic variety  is  isomorphic  to
 the initial one. It follows that $R(G)$  is an invariant of
 the  finitely generated presentation of the group
 $G$.  In the case of a free group of  rank  $n$, here denoted by  $ F_n$,
 the algebraic variety
 $R(F_n)$ is simply $(SL_2 \bold C)^n$,
  and is consequently an irreducible
 algebraic variety of dimension $3n$. Note that in this
 paper the term `algebraic variety' may be applied to both reducible and
 irreducible varieties.
 \medskip
 Let $V \subset \bold C^n$ be an algebraic variety (not necessarily
 reducible) and
 $V= V_1 \cup V_2 \dots \cup V_n$ be its unique decomposition into
 maximal irreducible components (see \cite {MD}).
 Then the number of maximal
 irreducible
 components of a given dimension is an invariant  of  $V$ and thus, if
 $V=R(G)$,  an  invariant  of  the  finitely  generated  group   $G$.
  \medskip
 The groups $G_{pt}$ are known to have $Dim(R(G))=4$, where $p,t$ are
 integers greater than one (see \cite {LS}, or the proof of Theorem A).
 It is the object of this communication
 to introduce this invariant and to
 produce a formula giving the
 number of irreducible 4-dimensional components of $R(G)$, where
 $G$ is isomorphic to one of the groups
 $$G_{pt}=\langle a, b; a^p= b^t \rangle. $$
 \medskip

 \proclaim  {Theorem  A}  Let  $C_4$  be the  number of
 irreducible  4 dimensional components of $R(G_{pt})$, where $p$ and $t$
 are integers greater than one. Then
 \roster
 \item "{a)}"  $C_4=$ $\frac {(p-2)(t-2)+pt} {4}$  if  both  $p,t$ are
 even.
 \item "{b)}"  $C_4=$ $\frac {(p-1)(t-1)} {2}$  if either  $p$ or $t$
 is odd.
 \endroster
 \endproclaim
 \medskip
 \flushpar
 That the class of groups $G_{pt}$ in Theorem 3.1  consists
 of  non-isomorphic  groups for different pairs
 $p,t \ge 2$  ({\it  up  to permutation of $p$ and $t$
 \/}) is a well known result  of  O.  Schreier  (1924),  \cite
 {SO}.

 \medskip
 Incidentally, a somewhat similar  but  different  study,  and  one  the
 author was unaware of when the main theorem of this note was obtained,  is
 counting
 connected components  in  the  space  of  representations  of  a  finitely
 generated group in a Lie group. For example, William Goldman in \cite {GW}
 develops a formula for counting the number of connected components in the
 space of representations of the fundamental group of an oriented surface
 in the $n$-fold covering group of $PSL(2,\bold R)$
 in terms of $n$ and the genus
 of the surface.
 \medskip
 Theorem A is of interest in its own right. However,  one  of  its
 consequences is rather curious and deserves mention  as
 it may suggest  perhaps   deeper  connections  between  a  group $G$
  and  the
 invariants of $R(G)$.  A well known fact   is   that
 a torus  knot  has the fundamental  group which   is  isomorphic  to  some
 $G_{pt}$, where $p$ and $t$ are relatively  prime   (see  \cite  {BZ}).  A
 direct consequence of Theorem A is that if $K$ is a  torus knot  then  its
 genus   equals   the   number   of   four   dimensional    components   in
 $R(G_{pt})$. The genus of a torus knot can also be defined  as
 half the rank of its commutator subgroup; consult \cite {BZ}.
 \medskip
 \head  Proofs \endhead
 \medskip
 \remark {Notation}
 Given a positive integer $p$ and $M \in SL_2 \bold C$, denote by
 $\Omega(p,M)$ the set
 $\{A \mid A\in SL_2 \bold C, A^p=M\}$. By $Tr(A)$ will be meant  the
 trace of the matrix $A$.
 Finally,
 denote by $I$ the $2 \times 2 $ identity matrix.
 \endremark
 \medskip

 \remark {Observation 1}
 $$R(G_{pt})=\{ (m_1,m_2) \vert m_1 \in SL_2 \bold C, \, m_2 \in
 \Omega(t,m_1^p) \}.$$
 \endremark
 \medskip
 The following proper sub-variety of $R(G_{pt})$ will play an important role
 in the proof of Theorem A.
 Let
 $$S=\{ (m_1, m_2 )\vert m_1^p=\pm I \}. \tag 1.0 $$
 Clearly
 $$S= S_+ \cup S_-,$$
 where
 $$S_+=\{ (m_1, m_2) \vert m_1^p= I \}, \,\,\,\,\,\,
 S_-=\{ (m_1, m_2 )\vert m_1^p= -I \}.$$
\medskip
\proclaim {Theorem B}
 $Dim((R(G_{pt})-S)) = 3$.
 \endproclaim
 \medskip
 \demo  {Proof}
 $\{ R(G_{pt})-S \}$ maps onto a quasi-affine variety $Q$ of $R(F_1)$
 via the map from $R(G_{pt})$ to $R(F_1)$  given by
 $ \phi(m_1, m_2)=m_1.$  Some care must be exercised in defining $Q$.
 If $p$ or $t$ is odd,  it can be assumed  without
 any loss of generality that $t$ is the odd one.
  This is done since if $t$ is even and $p$ odd, and  if
 $\,\,m_1^p \in \, Orb \,
 \left(\smallmatrix  -1 & 1 \\ 0 & -1 \endsmallmatrix\right) $,
 then
 $\Omega(t,m_1^p)= \emptyset$, where by $Orb$ is meant the orbit
 under conjugation by elements of  $SL_2 \bold C$; see
 \cite {LS}.
 Note that if $p$ is even and $m$ in $SL_2 \bold C$, then $m^p$ is never in
$Orb\,\left(\smallmatrix  -1 & 1 \\ 0 & -1 \endsmallmatrix\right)$.
 This said, now proceed by letting
 $$Q= R(F_1) - \phi(S), $$
 where $S$ is as in (1.0).
\medskip
 A few comments are in line regarding fibres of the map $\phi$.
 If $m \in
 \phi (S)$ and $m^p=-I$, then $Dim (\phi^{-1}(m))=2$. This follows by Lemma
 1.6 below.  Similarly, if $m \in
 \phi (S)$ and $m^p=I$, then $Dim (\phi^{-1}(m))=2$ whenever $p\ge 3$.
 This is
 a consequence of Lemma 1.6 below. Any $m \in Q$
 has the property that $Dim (\phi^{-1}(m))=0$.  This is a result of
 the fact that if $\,\,m\in \{SL_2 \bold C- \{\pm I \}\}$, and $t$ an
 integer, then the
 equation
 $\,\,x^t=m$ has at
 most  a  finite  number of solutions in  $SL_2 \bold C$; see  \cite  {LS},
 or \cite {GF}.
 As a consequence (by
 elementary arguments involving fibres of regular maps with a dense image
 between varieties (see \cite {MD}, or \cite {LS} ),
 since
 the  fibres  of  $\phi$  over  the  quasi-affine  variety  $Q$  are   zero
 dimensional, and $(R(G_{pt})-S)$ maps onto the  quasi-affine  variety  $Q$,
 one obtains that
 $$Dim(Q)= Dim((R(G_{pt})-S)).$$
 Notice that $Q$ is a quasi-affine variety in the irreducible variety
 $SL_2 \bold C$, and
 thus $Dim(Q)=Dim (SL_2 \bold C)=3$. This concludes the proof of Theorem B.
 \enddemo

 \medskip
 It follows then that if any four dimensional components are present
 in $R(G_{pt})$, then they are in the sub-variety  $S$.
 \medskip

 Next we propose to count all the four dimensional irreducible components
 in $S$, for different values of $p,t$. We will find necessary  the
 following elementary lemmas stated here mostly for the convenience of the
 reader;
 their demonstrations
 can be obtained using  naive
 facts from, for example,  \cite {GF} and \cite {MD}, or they can be
 found directly
 in \cite {LS}.
 \medskip
 \proclaim { Lemma 1.4 }
 Let $A \, \in \, SL_2 \bold C$
 be any matrix of a given trace $b\ne\pm 2$.
 Then any
 matrix $B$ in $SL_2 \bold C$
 having trace $b$ is similar to $A$.
 \endproclaim
 \medskip
 \proclaim {Lemma 1.5}
 Let  $A \in\,SL_2 \bold C$  be  of
  trace  $b\ne\pm2$.  Then  the  orbit  of
 $A$
 under  $SL_2 \bold C$ conjugation is an irreducible
 affine algebraic variety
  of  dimension 2.
 \endproclaim
 \medskip
 \proclaim {Lemma  1.6}
 \roster
 \item  "{i)}"  If  $p=2$, then
 $Dim\,\Omega(p,I)=0$,  and   $\Omega(p,I)$ is
 reducible.
 \item "{ii)}" If $p>2$, then
 $Dim\,\Omega(p,I)=2$, and $\Omega(p,I)$ is
 reducible.
 \item "{ iii)}" If $p\ge2$, $Dim\, \Omega(p,-I)=2$, and for
 for $p > 2$, $\Omega(p,-I)$ is reducible.
 \endroster
 \endproclaim
 \medskip

  From Observation 1 together with (1.0), it is easy to deduce that
 $$ S_+ =\Omega(p,I) \times \Omega(t,I). \tag 1.1 $$
 $$ S_- =\Omega(p,-I) \times \Omega(t,-I). \tag 1.2 $$
 \medskip
 Using Lemmas 1.4 and 1.5 it is possible to deduce that the number of
 two dimensional irreducible components in $\Omega (p,I)$ is given by:
 $$ \frac {p-1} {2},\,\, \text {if} \,\, p \,\, \text {is odd,} \tag 1.3 $$
 $$ \frac {p-2} {2},\,\, \text {if} \,\, p \,\, \text {is even.} \tag 1.4 $$
 In a similar fashion one can be deduced that the number of two
 dimensional irreducible components in $\Omega (p,-I)$ is given by:
 $$ \frac {p-1} {2},\,\, \text {if} \,\, p \,\, \text {is odd,} \tag 1.5 $$
 $$ \frac {p} {2},\,\, \text {if} \,\, p \,\, \text {is even.} \tag 1.6 $$
 \medskip
 Note, having made these deductions, it is clear that the $R(G_{pt})$ are
 four  dimensional reducible varieties.
 \medskip
 \demo { Proof of Theorem A }
 If both $p$ and $t$ are even then using (1.1) and (1.2) together
 with (1.4) and (1.6)  the number of four dimensional components
 in $R(G_{pt})$ is
 given by:
 $$ (\frac {p-2} {2})(\frac {t-2} {2}) +
 (\frac {p} {2})(\frac {t} {2}).$$
 \medskip
 If both $p$ and $t$ are odd then the number of four dimensional components
 is given by:
$$ (\frac {p-1} {2})(\frac {t-1} {2}) +
 (\frac {p-1} {2})(\frac {t-1} {2}).$$
 \medskip
 In the case that $p$ is odd and $t$ is even then the number of four
 dimensional components is  given by:
 $$ (\frac {p-1} {2})(\frac {t-2} {2}) +
 (\frac {p-1} {2})(\frac {t} {2}).$$
 The proof of Theorem A is now complete  since elementary algebra
 yields the quantities promised in the statement of Theorem A.
 \enddemo
 
 \enddocument